\documentclass[11pt]{article}
\usepackage[mathscr]{eucal}
\usepackage{mathrsfs}
\usepackage{amssymb}
\usepackage{calligra}
\usepackage{fontenc}
\usepackage{amsbsy} 
\usepackage{graphicx}
\graphicspath{%
    {converted_graphics/}
    {/}
}
\begin{document}

\newcommand{\ts}{{\tilde{\sf s}}}
\newcommand{\sfv}{{\sf v}}
\newcommand{\sfw}{{\sf w}}
\newcommand{\simge}{\ba{cc}\vspace*{-2.4mm}>\\ \sim\ea }
\newcommand{\simle}{\ba{cc}\vspace*{-2.4mm}<\\ \sim\ea }
\newcommand{\Cdot}{\!\cdot\!}
\newcommand{\sq}{{$\sqcap\!\!\!\!\sqcup$}}
\newcommand{\Eu}{{\rm I\,\!\! E}}
\newcommand{\Io}{\Int{\Omega}{}}
\newcommand{\Id}{\Int{\cald}{}}
\newcommand{\Div}{\mbox{\rm div}\,}
\newcommand{\tr}{\mbox{\rm tr}\,}
\newcommand{\grad}{\mbox{\rm grad}\,}
\newcommand{\supp}{\mbox{\rm supp}\,}
\newcommand{\curl}{\mbox{\rm curl}\,}
\newcommand{\Ido}{\Int{\partial\Omega}{}}
\newcommand{\IdS}{\Int{\Sigma}{}}
\newcommand{\Oint}[2]{{\displaystyle \oint_{#1}^{#2}}}
\newcommand{\Int}[2]{{\displaystyle \int_{ #1}^{ #2}}}
\newcommand{\Lim}[1]{{\displaystyle \lim_{ #1}}}
\newcommand{\Limsup}[1]{{\displaystyle \limsup_{\footnotesize #1}}}
\newcommand{\Liminf}[1]{{\displaystyle \liminf_{\footnotesize #1}}}
\newcommand{\Sup}[1]{{\displaystyle \sup_{#1}}}
\newcommand{\Inf}[1]{{\displaystyle \inf_{#1}}}
\newcommand{\Max}[1]{{\displaystyle \max_{#1}}}
\newcommand{\Min}[1]{{\displaystyle \min_{#1}}}
\newcommand{\Sum}[2]{{\displaystyle \sum_{#1}^{#2}}}
\newcommand{\Prod}[2]{{\displaystyle \prod_{#1}^{#2}}}
\newcommand{\BCup}[2]{{\displaystyle \bigcup_{#1}^{#2}}}
\newcommand{\BCap}[2]{{\displaystyle \bigcap_{#1}^{#2}}}
\newcommand{\Frac}[2]{\displaystyle{\frac{\displaystyle{#1}}{\displaystyle{#2}}}}
\newcommand{\norm}[1]{\left\|{#1}\right\|}
\newcommand{\Norm}[1]{\langle\langle{#1}\rangle\rangle_q}
\newcommand{\No}[1]{\langle\!\langle{#1}\rangle\!\rangle}
\newcommand{\NO}[1]{{\langle{#1}\rangle}_{\lambda,q}}
\newcommand{\beea}{\begin{eqnarray}}
\newcommand{\eeea}{\end{eqnarray}}
\newcommand{\ms}{\medskip\smallskip}
\newcommand{\bs}{\bigskip}
\newcommand{\ps}{\par\smallskip}
\newcommand{\bfe}{{\mbox{\boldmath $e$}} }
\newcommand{\pni}{{\par\noindent}}
\newcommand{\bfq}{{\mbox{\boldmath $q$}} }
\newcommand{\bfz}{{\mbox{\boldmath $z$}} }
\newcommand{\0}{{\mbox{\boldmath $0$}} }
\newcommand{\LE}{\!\!\!&\le&\!\!\!}
\newcommand{\BL}[1]{{\par\smallskip{\bf Lemma #1.}}}
\newcommand{\BT}[1]{{\par\smallskip{\bf Theorem #1.}}}
\newcommand{\Ln}{[\!|}
\newcommand{\Rn}{|\!]}
\newcommand{\n}[1]{{\Ln{#1}\Rn}} 
\newcommand{\nq}[1]{{\Ln{#1}\Rn}_{q}} 
\newcommand{\nqr}[1]{{\Ln{#1}\Rn}_{q,r}} 
\newcommand{\Nq}[1]{{\langle{#1}\rangle}_{q}} 
\newcommand{\Nql}[1]{{\langle{#1}\rangle}_{\lambda,q}} 
\newcommand{\Nqr}[1]{{\langle{#1}\rangle}_{q,r}}
\newcommand{\N}[1]{{|\!\!|\!\!|\,{#1}\,|\!|\!\!|_2}}
\newcommand{\EA}[2]{$$#1$$%
\vspace{-6.mm}
\begin{equation}
\end{equation}
\vspace{-6.mm}
$$
#2
\setlength{\belowdisplayskip}{3mm}
\setlength{\belowdisplayshortskip}{3mm}
$$
}
\newcommand{\A}[2]{$$#1$$%
\vspace{-4.mm}
$$
#2
\setlength{\belowdisplayskip}{3mm}
\setlength{\belowdisplayshortskip}{3mm}
$$
}
\newcommand{\BF}{\begin{footnotesize}}
\newcommand{\EF}{\end{footnotesize}}
\setlength{\jot}{.15in}
\newcommand{\pde}[2]{{\displaystyle \frac{\mbox{$\partial #1$}}{\mbox{$\partial #2$}}}}
\newcommand{\ode}[2]{{\displaystyle \frac{\mbox{$d #1$}}{\mbox{$d #2$}}}}
\newcommand{\f}[2]{\frac{\mbox{$#1$}}{\mbox{$ #2$}}}
\newcommand{\bi}{\begin{itemize}}
\newcommand{\ei}{\end{itemize}}
\newcommand{\ed}{\end{document}}
\newcommand{\be}{\begin{equation}}
\newcommand{\ba}{\begin{array}}
\newcommand{\ea}{\end{array}}
\newcommand{\ee}{\end{equation}}
\newcommand{\eeq}[1]{\label{eq:#1}\end{equation}}
\newcommand{\real}{{\mathbb R}}
\newcommand{\compl}{{\mathbb C}}
\def\Id{\mbox{\boldmath $1$}}
\def\zero{\mbox{\boldmath $0$}}
\newcommand{\PP}{{\rm I\!\!\,P}}
\newcommand{\nat}{{\mathbb N}}
\newcommand{\bfpsi}{\mbox{\boldmath $\psi$}}
\newcommand{\bfchi}{\mbox{\boldmath $\chi$}}
\newcommand{\bfomega}{\mbox{\boldmath $\omega$}}
\newcommand{\bfome}{\mbox{\boldmath $\varpi$}}
\newcommand{\bfvaromega}{\mbox{\boldmath $\varpi$}}
\newcommand{\bfOmega}{\mbox{\boldmath $\Omega$}}
\newcommand{\bfTheta}{\mbox{\boldmath $\Theta$}}
\newcommand{\bfxi}{\mbox{\boldmath $\xi$}}
\newcommand{\bfmu}{\mbox{\boldmath $\mu$}}
\newcommand{\bfx}{\mbox{\boldmath $x$}}
\newcommand{\bfy}{\mbox{\boldmath $y$}}
\newcommand{\bfPsi}{\mbox{\boldmath $\Psi$}}
\newcommand{\bfphi}{\mbox{\boldmath $\varphi$}}
\newcommand{\bfhi}{\mbox{\boldmath $\phi$}}
\newcommand{\bfPhi}{\mbox{\boldmath $\Phi$}}
\newcommand{\bfv}{{\mbox{\boldmath $v$}} }
\newcommand{\bfu}{{\mbox{\boldmath $u$}} }
\newcommand{\bfsf}{{\mbox{\footnotesize\boldmath $s$}} }
\newcommand{\bfuf}{{\mbox{\footnotesize\boldmath $u$}} }
\newcommand{\bfw}{{\mbox{\boldmath $w$}} }
\newcommand{\bff}{{\mbox{\boldmath $f$}} }
\newcommand{\bfa}{{\mbox{\boldmath $a$}} }
\newcommand{\bfi}{{\mbox{\boldmath $i$}} }
\newcommand{\bfj}{{\mbox{\boldmath $j$}} }
\newcommand{\bfc}{{\mbox{\boldmath $c$}} }
\newcommand{\bfo}{{\mbox{\boldmath $o$}} }
\newcommand{\bfp}{{\mbox{\boldmath $p$}} }
\newcommand{\bfkp}{{\mbox{\footnotesize{\boldmath $k$}}} }
\newcommand{\bfka}{{\mbox{\footnotesize{\boldmath $k^*$}}} }
\newcommand{\bft}{{\mbox{\boldmath $t$}} }
\newcommand{\bfd}{{\mbox{\boldmath $d$}} }
\newcommand{\bfl}{{\mbox{\boldmath $l$}} }
\newcommand{\bfr}{{\mbox{\boldmath $r$}} }
\newcommand{\bfk}{{\mbox{\boldmath $k$}} }
\newcommand{\bfA}{{\mbox{\boldmath $A$}} }
\newcommand{\bfS}{{\mbox{\boldmath $S$}} }
\newcommand{\bfO}{{\mbox{\boldmath $O$}} }
\newcommand{\bfM}{{\mbox{\boldmath $M$}} }
\newcommand{\bfP}{{\mbox{\boldmath $P$}} }
\newcommand{\bfB}{{\mbox{\boldmath $B$}} }
\newcommand{\bfR}{{\mbox{\boldmath $R$}} }
\newcommand{\bfC}{{\mbox{\boldmath $C$}} }
\newcommand{\bfD}{{\mbox{\boldmath $D$}} }
\newcommand{\bfQ}{{\mbox{\boldmath $Q$}} }
\newcommand{\bfZ}{{\mbox{\boldmath $Z$}} }
\newcommand{\bfG}{{\mbox{\boldmath $G$}} }
\newcommand{\bfE}{{\mbox{\boldmath $E$}} }
\newcommand{\bfX}{{\mbox{\boldmath $X$}} }
\newcommand{\bfY}{{\mbox{\boldmath $Y$}} }
\newcommand{\bfH}{{\mbox{\boldmath $H$}} }
\newcommand{\bfI}{{\mbox{\boldmath $I$}} }
\newcommand{\bfJ}{{\mbox{\boldmath $J$}} }
\newcommand{\bfN}{{\mbox{\boldmath $N$}} }
\newcommand{\bfh}{{\mbox{\boldmath $h$}} }
\newcommand{\bfm}{{\mbox{\boldmath $m$}} }
\newcommand{\bfone}{{\mbox{\boldmath $1$}} }
\newcommand{\hs}{{\rm I}\!\!\,{\rm R}^3_+}
\newcommand{\cala}{{\cal A}}
\newcommand{\calb}{{\cal B}}
\newcommand{\calc}{{\cal C}}
\newcommand{\cald}{{\cal D}}
\newcommand{\cale}{{\cal E}}
\newcommand{\calf}{{\cal F}}
\newcommand{\calg}{{\cal G}}
\newcommand{\calh}{{\cal H}}
\newcommand{\cali}{{\cal I}}
\newcommand{\calj}{{\cal J}}
\newcommand{\calk}{{\cal K}}
\newcommand{\call}{{\cal L}}
\newcommand{\calm}{{\cal M}}
\newcommand{\caln}{{\cal N}}
\newcommand{\calo}{{\cal O}}
\newcommand{\calp}{{\cal P}}
\newcommand{\calq}{{\cal Q}}
\newcommand{\calr}{{\cal R}}
\newcommand{\cals}{{\cal S}}
\newcommand{\calt}{{\cal T}}
\newcommand{\calu}{{\cal U}}
\newcommand{\calv}{{\cal V}}
\newcommand{\calx}{{\cal X}}
\newcommand{\caly}{{\cal Y}}
\newcommand{\calw}{{\cal W}}
\newcommand{\calz}{{\cal Z}}
\newcommand{\bfsigma}{\mbox{\boldmath $\sigma$}}
\newcommand{\bfSigma}{\mbox{\boldmath $\Sigma$}}
\newcommand{\bftau}{\mbox{\boldmath $\tau$}}
\newcommand{\bfeta}{\mbox{\boldmath $\eta$}}
\newcommand{\bfT}{{\mbox{\boldmath $T$}} }
\newcommand{\bfV}{{\mbox{\boldmath $V$}} }
\newcommand{\bfU}{{\mbox{\boldmath $U$}} }
\newcommand{\bfW}{{\mbox{\boldmath $W$}} }
\newcommand{\bfF}{{\mbox{\boldmath $F$}} }
\newcommand{\bfK}{{\mbox{\boldmath $K$}} }
\newcommand{\bfL}{{\mbox{\boldmath $L$}} }
\newcommand{\bfb}{{\mbox{\boldmath $b$}} }
\newcommand{\bfg}{{\mbox{\boldmath $g$}} }
\newcommand{\bfn}{{\mbox{\boldmath $n$}} }
\newcommand{\bfs}{{\mbox{\boldmath $s$}} }
\newcommand{\cf}{{\it cf.} }
\newcommand{\io}{\int_\Omega}
\newcommand{\1}{\item[({\it i})]}
\newcommand{\2}{\item[({\it ii})]}
\newcommand{\3}{\item[({\it iii})]}
\newcommand{\4}{\item[({\it iv})]}
\newcommand{\5}{\item[({\it v})]}
\newcommand{\6}{\item[({\it vi})]}
\newcommand{\7}{\item[({\it vii})]}
\newcommand{\8}{\item[({\it viii})]}
\newcommand{\9}{\item[({\it xi})]}
\newcommand{\ido}{\int_{\partial\Omega}}
\newcommand{\half}{\mbox{$\frac{1}{2}$}}
\def\parallel{\|}
\def\mid{|}
\def\Bbb R{\real}
\def\hat{\widehat}
\def\tilde{\widetilde}
\def\bar{\overline}
\newcommand{\threehalves}{3\over 2}
\newcommand{\bfPi}{\mbox{\boldmath $\Pi$}}
\newcommand{\bfXi}{\mbox{\boldmath $\Xi$}}
\newcommand{\bfalpha}{\mbox{\boldmath $\alpha$}}
\newcommand{\bfbeta}{\mbox{\boldmath $\beta$}}
\newcommand{\bfgamma}{\mbox{\boldmath $\gamma$}}
\newcommand{\bfdelta}{\mbox{\boldmath $\delta$}}
\newcommand{\bfzeta}{\mbox{\boldmath $\zeta$}}
\newcommand{\bfUpsilon}{\mbox{\boldmath $\Upsilon$}}
\newcommand{\bfGamma}{\mbox{\boldmath $\Gamma$}}
\newcommand{\bfcala}{\mbox{\boldmath ${\cal A}$}}
\newcommand{\bfcalm}{\mbox{\boldmath ${\cal M}$}}
\newcommand{\bfcaln}{\mbox{\boldmath ${\cal N}$}}
\newcommand{\bfcalq}{\mbox{\boldmath ${\cal Q}$}}
\newcommand{\bfcalb}{\mbox{\boldmath ${\cal B}$}}
\newcommand{\bfcalc}{\mbox{\boldmath ${\cal C}$}}
\newcommand{\bfcali}{\mbox{\boldmath ${\cal I}$}}
\newcommand{\bfcalg}{\mbox{\boldmath ${\cal G}$}}
\newcommand{\bfcalh}{\mbox{\boldmath ${\cal H}$}}
\newcommand{\bfcalk}{\mbox{\boldmath ${\cal K}$}}
\newcommand{\bfcalt}{\mbox{\boldmath ${\cal T}$}}
\newcommand{\bfcalx}{\mbox{\boldmath ${\cal X}$}}
\newcommand{\bfcall}{\mbox{\boldmath ${\cal L}$}}
\newcommand{\bfcalf}{\mbox{\boldmath ${\cal F}$}}
\newcommand{\bfcalr}{\mbox{\boldmath ${\cal R}$}}
\newcommand{\bfcals}{\mbox{\boldmath ${\cal S}$}}
\newcommand{\bfcalw}{\mbox{\boldmath ${\cal W}$}}
\newcommand{\bfcalu}{\mbox{\boldmath ${\cal U}$}}
\newcommand{\bfcalv}{\mbox{\boldmath ${\cal V}$}}
\newcommand{\bfcalz}{\mbox{\boldmath ${\cal Z}$}}
\pagenumbering{roman}
\newcommand{\art}[6]{{\I[{\sc #1,}] {#2}, {\it #3}, {\bf #4}, {#5} {[#6]}}}
\newcommand{\ED}{\end{description}}
\newcommand{\I}{\item }
\newcommand{\ra}{\rm a}
\newcommand{\rb}{\rm b}
\newcommand{\rc}{\rm c}
\newcommand{\Hsp}{{\rm I}\!\!\,{\rm R}^n_+}
\newcommand{\Hsn}{{\rm I}\!\!\,{\rm R}^n_-}
\newcommand{\po}[1]{\mbox{$\displaystyle \frac{\mbox{$\partial #1$}}
{\mbox{$\partial x_{1}$}}$}}
\newcommand{\PO}[1]{\mbox{$\displaystyle \frac{\mbox{$\partial #1$}}
{\mbox{$\partial y_{1}$}}$}}
\newcommand{\OP}{\left(\Delta+2\lambda\PO{}\right)}
\newcommand{\op}{\left(\Delta+2\lambda\po{}\right)}
\newcommand{\ft}[1]{
\Frac{1}{(2\pi)^{n/2}}\Int{{\Bbb R}^{n}}{}e^{i{\bf x}\cdot \bfxi}
#1(\xi)d\xi}
\newcommand{\Ft}[1]{
\Frac{1}{2\pi}\Int{{\Bbb R}^{2}}{}e^{i{x}\cdot \xi}
#1(\xi)d\xi}
\newcommand{\Z}{\item[({\it a})]}
\newcommand{\B}{\item[({\it b})]}
\newcommand{\C}{\item[({\it c})]}
\newcommand{\D}{\item[({\it d})]}
\newcommand{\E}{\item[({\it e})]}
\newcommand{\G}{\item[({\it g})]}
\newcommand{\Š}{\`e}
\newcommand{\…}{\`a}
\newcommand{\•}{\`o}
\newcommand{\—}{\`u}
\newcommand{\}{\`{\i}}
\def\tag{\renewcommand{\theequation}}
\newcommand{\Footnote}{~\footnote}
\newcommand{\ie}{{\it i.e.}}
\newcommand{\dist}{\mbox{\rm dist\,}}
\newcommand{\const}{\mbox{\rm const}}
\newcommand{\trace}{\mbox{\rm trace}}
\newcommand{\Bo}{\par\hfill{$\Box$}\par\noindent}
\newcommand{\Nor}[1]{\langle{#1}\rangle_q}
\newcommand{\vs}{\vspace*{.5cm}\par\noindent}
\newcommand{\Vs}{\vspace*{.6cm}\par\noindent}
\newcommand{\Vvs}{\vspace*{.7cm}\par\noindent}
\newcommand{\VVs}{\vspace*{.8cm}\par\noindent}
\newtheorem{definition}{Definition}[section]
\newcommand{\Bd}{\begin{definition}\begin{rm}}
\newcommand{\Ed}{\end{rm}\end{definition}}
\newtheorem{remark}{Remark}[section]
\newcommand{\Br}{\begin{remark}\begin{rm}}
\newcommand{\Er}{\end{rm}\end{remark}}
\newtheorem{proposition}{Proposition}[section]
\newcommand{\Bp}{\begin{proposition}\begin{sl}}
\newcommand{\EP}[1]{\end{sl}\label{proposition:#1}\end{proposition}}
\newcommand{\propref}[1]{{\rm Proposition \ref{proposition:#1}}}
\newcommand{\Bt}{\begin{theorem}\begin{sl}}
\newcommand{\Et}{\end{sl}\end{theorem}}
\newcommand{\Bl}{\begin{lemma}\begin{sl}}
\newcommand{\El}{\end{sl}\end{lemma}}
\newtheorem{theorem}{Theorem}[section]
\newtheorem{lemma}{Lemma}[section]
\newtheorem{corollary}{Corollary}[section]
\newcommand{\eqref}[1]{{\rm (\ref{eq:#1})}}
\newcommand{\Bc}{\begin{corollary}\begin{sl}}
\newcommand{\Ec}{\end{sl}\end{corollary}}
\newcommand{\ET}[1]{\end{sl}\label{theorem:#1}\end{theorem}}
\newcommand{\EDD}[1]{\end{rm}\label{definition:#1}\end{definition}}
\newcommand{\EL}[1]{\end{sl}\label{lemma:#1}\end{lemma}}
\newcommand{\theoref}[1]{{\rm Theorem \ref{theorem:#1}}}
\newcommand{\ER}[1]{\end{rm}\label{remark:#1}\end{remark}}
\newcommand{\EC}[1]{\end{sl}\label{corollary:#1}\end{corollary}}
\newcommand{\remref}[1]{{\rm Remark \ref{remark:#1}}}
\newcommand{\cororef}[1]{{\rm Corollary \ref{corollary:#1}}}
\newcommand{\lemmref}[1]{{\rm Lemma \ref{lemma:#1}}}
\newcommand{\essup}[1]{{\rm ess}\,{{\displaystyle \sup_{\hspace*{-5mm}{#1}}}}}

\renewcommand{\i}{{\rm i}}

\pagenumbering{arabic}
\newcommand{\QED}{{\par\hfill$\square$\par}}
\renewcommand{\thefootnote}{(\arabic{footnote})}
\title{On the Energy Equality for Distributional Solutions to Navier-Stokes Equations} 
\author{ Giovanni P. Galdi 
\thanks{Department of Mechanical Engineering and Materials Science, University of Pittsburgh, PA 15261. 
Work  partially supported by NSF DMS Grant-1614011.}}
\date{}
\maketitle
\begin{abstract}  A classical result of J.-L. Lions asserts that if  a solution to the Navier-Stokes equations is such that: (i) it is in the Leray-Hopf class, and (ii)  belongs to $L^4(0,T;L^4)$, then it must satisfy the energy equality in the time interval $[0,T]$. In this note we show that assumption (i) is not necessary. 
 \end{abstract}

\renewcommand{\theequation}{\arabic{section}.\arabic{equation}}
\setcounter{section}{0}
\section{Introduction} 
Consider the three-dimensional Cauchy problem\footnote{Generalizations to other domains are discussed in \remref{1}. We also assume, for simplicity, zero body force and,  without loss of generality, take the kinematic viscosity coefficient to be 1.} for the Navier-Stokes equations: 
\be\ba{cc}\smallskip\left.\ba{ll}\smallskip
\partial_tv+v\cdot\nabla v=\Delta v-\nabla p\\
\Div v=0\ea\right\}\ \ \mbox{in}\,\ \real^3\times (0,\infty)
\\
v(x,0)=v_0(x)\ \ x\in\real^3\,,
\ea
\eeq{1.1}
where $v:\real^3\times  
[0,\infty)\to v(x,t)\in\real^3$ is the flow velocity field and $p$ the associated pressure field. \par
It is  well known since the work of Leray \cite{Leray} and Hopf \cite{Hopf}, that for any $v_0\in L^2_\sigma(\real^3)$ one can construct a global weak solutions  to \eqref{1.1}, namely, a function $v$ that, for each $T>0$, is in the class\footnote{$L^2_\sigma(\real^3)$ is the subspace of $L^2(\real^3)$ of divergence-free vector functions. Other notations are standard, like $H^{m,q}$,  for Sobolev spaces, with corresponding norm $\|\cdot\|_{m,q}$, $L^r(I;X)$, $I$  real interval, $X$ Banach space, for Bochner spaces, etc. }
\be
v\in L^\infty(0,T;L^2_\sigma(\real^3))\cap L^2(0,T;H^{1}(\real^3))
\eeq{1.2}
and solves \eqref{1.1} in a distributional sense. In addition, such a $v$ satisfies the so-called energy inequality:\footnote{Actually, $v$ obeys the {\em strong energy inequality} \cite[\S\S 27--28]{Leray}, but this is irrelevant to the aim of our paper.} 
\be
\|v(t)\|_2^2+2\int_0^t\|\nabla v(\tau)\|_2^2\le\|v_0\|_2^2\,,\ \ \mbox{all $t\ge 0$}\,,
\eeq{1.3}
where $\|\cdot\|_q$,  denotes the $L^q(\real^3)$ -norm. The inequality sign in this relation is another indication of poor regularity of a weak solution. Actually, all sufficiently regular solutions to \eqref{1.1} satisfy \eqref{1.3} with the equality sign ({\em energy equality}), which provides the precise energy balance for the given flow. As a matter of fact, it is still an outstanding open question whether there exist  {\em global} solutions to \eqref{1.1} satisfying  the energy equality for arbitrary $v_0\in L^2_\sigma$. \par
In this regard, a famous result of J.-L. Lions     \cite{Lions} states that if $v$, in {\em addition} to be in the class \eqref{1.2},   satisfies also
\be
v\in L^4(0,T;L^4(\real^3))\,,
\eeq{H}
then necessarily $v$ obeys the energy equality throughout the interval $[0,T]$.\footnote{For further sufficient conditions other than \eqref{H}, see \cite{ChesFri} and the references therein.} An important corollary to this finding is that the map $t\mapsto v(t)$ is continuous with values in $L^2_\sigma$. Notice that the two classes \eqref{1.2} and \eqref{H} are not comparable, in the sense that a generic function in \eqref{1.2} need not satisfy \eqref{H} and vice versa.
\par
Objective of this note is to prove that, actually,  for the validity of Lions result the requirement  \eqref{1.2} is entirely redundant:  it is {\em just enough} that $v$ satisfies \eqref{H},  along with the (necessary) condition $v_0\in L^2_\sigma(\real^3)$. More precisely, setting
$$
\cald_T:=\{\varphi\in C_0^\infty(\real^3\times [0,T)):\ \Div\varphi=0\}
$$
we will show the following. 
\Bt Let $v\in L^2_{\rm loc}(\real^3\times(0,T))$ be such that  
\be\ba{ll}\medskip
\Int0T\Int{\real^3}{}v\cdot(\partial_t\varphi+\Delta\varphi+v\cdot\nabla\varphi)=-\Int{\real^3}{}v_0\cdot\varphi(0)\\
\Int0T\Int{\real^3}{}v\cdot\phi=0
\ea
\eeq{1.5}
for some $v_0\in L^2_\sigma(\real^3)$ and all $\varphi\in \cald_T$, $\phi\in C_0^\infty(\real^3\times(0,T))$.
Then, if $v$ satisfies
 \eqref{H}, necessarily $v$ is in the class \eqref{1.2} and thus obeys the energy equality. 
\ET{1}
\par
The proof of this result, which we give in the following section, is surprisingly elementary, and is based on a mollifying procedure coupled with a simple duality argument.
\setcounter{equation}{0}
\section{Proof of Theorem 1.1}
For a given $g:\real^3\times (0,T)\mapsto\real^3$ we set
$
\tilde{g}=g(\cdot,T-t)\,.
$
If $g$ is locally integrable, by $g^{(\eta)}$, and $g_{(\eta)}$, $\eta$ positive and sufficiently small, we denote  space and space-time mollifiers of $g$:  
$$
g^{(\eta)}(x,\cdot)=\Int{\real^3}{}k_\eta(x-y)g(y,\cdot)dy\,,\ \
g_{(\eta)}(x,t)=\Int0Tj_\eta(t-s)g^{(\eta)}(x,s)ds
\,,
$$
where 
$$
j_\eta(\tau):=\eta^{-1}j(\tau/\eta)\,,\ \ k_\eta(\xi):=\eta^{-1}k(\xi/\eta)\,,\ \ (\tau,\xi)\in \real\times\real^3\,, 
$$
with $j\in C_0^\infty(-1,1)$, and $k\in C_0^\infty(\real^3)$. 
We also write, as customary,
$
(f,g):=\int_{\real^3}f\cdot g\,,
$
and set, for simplicity,  $L^{p,q}:=L^p(0,T;L^q(\real^3))$. Finally, we define
$$
\mathscr W^{1,p}:=\big\{u\in L_{\rm 1oc}^1(\real^3\times (0,T)): u\in H^{1,p}(0,T;L^p(\real^3))\cap L^p(0,T;H^{2,p}(\real^3))\big\}\,.
$$
\par Before proving the theorem, we need a preparatory result.
\Bl Let $v$ and $v_0$ be as in \theoref{1} and $f\in C_0^\infty(\real^3\times(0,T))$. Then the Cauchy problem\be\ba{cc}\smallskip\left.\ba{ll}\smallskip
\partial_tw+\tilde{ v_{(\eta)}}\cdot\nabla w=\Delta w-\nabla {\sf p}+\tilde{f}\\
\Div w=0\ea\right\}\ \ \mbox{in}\, \real^3\times (0,\infty)
\\
w(x,0)=(v_0)^{(\eta)}(x)\ \ x\in\real^3\,,
\ea
\eeq{2.1}
has one (and only one) solution $(w_\eta,{\sf p}_\eta)$ such that\footnote{The continuous embedding in \eqref{2.2} is a classical interpolation result \cite{Lions1}.}
\be
w_\eta \in \mathscr W^{1,2}\subset C([0,T];H^1)\,,\ \ \nabla{\sf p}_\eta\in  L^{2,2}\,,
\eeq{2.2}
and satisfying the uniform bound
\be
\max_{t\in [0,T]}\|w_\eta(t)\|_2^2+\int_0^T\|\nabla w_\eta(t)\|_2^2\le \|v_0\|_2^2+c_2\int_0^T\|f(t)\|_{\frac65}^2\,.
\eeq{2.4}
In addition, if $v_0\equiv 0$, we have also $(w_\eta,\nabla{\sf p}_\eta)\in \mathscr W^{1,\frac43}\times L^{\frac43,\frac43}$.
\EL{2.2}
{\em Proof.} The existence of a solution in the class \eqref{2.2} is easily established by the ``invading domains" technique coupled with the classical Galerkin method \cite{Hey}. We will sketch a proof here. Let $B_R\subset\real^3$ be the ball of radius $R$ centered at the origin, and consider the following problem:
\be\ba{cc}\smallskip\left.\ba{ll}\smallskip
\partial_tw+\tilde{ v_{(\eta)}}\cdot\nabla w=\Delta w-\nabla {\sf p}+\tilde{f}\\
\Div w=0\ea\right\}\ \ \mbox{in}\, B_R\times (0,\infty)
\\
w(x,t)=0\ \ (x,t)\in\partial B_R\times (0,\infty)\,,\ \ w(x,0)=v_R(x)\ \ x\in B_R\,,
\ea
\eeq{2.5}
where $v_R\in H^1_0(B_R)\cap L^2_\sigma(B_R)$ satisfies\footnote{Since $v_0\in L^2_\sigma$, we have  $(v_0)^{(\eta)}\in H^{m,q}$, for all $m\ge 0$ and all $q\in [2,\infty]$.}
\be
\lim_{R\to\infty}\|v_R-(v_0)^{(\eta)}\|_{1,2}=0\,.
\eeq{2.6} 
The existence of  $v_R$ is known \cite[Theorem III.4.2]{GaBook}. If we dot-multiply both sides of \eqref{2.5}$_1$ by $w$, formally integrate by parts over $B_R$ and take into account \eqref{2.5}$_2$ and  $\Div\tilde{ v_{(\eta)}}=0$, we get
$$
\|w(t)\|_2^2+2\int_0^t\|\nabla w(\tau)\|_2^2=\|v_R\|_2^2+\int_0^t(\tilde f(\tau),w(\tau))\,
$$
The latter, in turn, by \eqref{2.6},  Sobolev inequality $\|w\|_6\le c_1\,\|\nabla w\|_2$, and Cauchy--Schwartz inequality furnishes
\be
\|w(t)\|_2^2+\int_0^t\|\nabla w(\tau)\|_2^2\le \|v_0\|_2^2+c_2\int_0^t\|f(\tau)\|_{\frac65}^2\,.
\eeq{2.7}
Furthermore, since
$
\tilde{ v_{(\eta)}}\in L^{\infty,\infty}\,,
$
if we dot-multiply both sides of \eqref{2.5}$_1$ a first time by  $P\Delta w$, with $P:L^2\mapsto L^2_\sigma$ Helmholtz projector, and then a second time by $\partial_tw$,  and then integrate by parts the resulting relations over $B_R$, we get
$$\ba{rl}\medskip
\half\ode{}t\|\nabla w(t)\|_2^2+\|P\Delta w\|_2^2=&\!\!\!\!\!
(\tilde{ v_{(\eta)}}\cdot\nabla w,P\Delta w)+(\tilde{f},w)\\
\le &\!\!\!\!\!
\left(\|\tilde{ v_{(\eta)}}\|_{L^{\infty,\infty}}\|\nabla w\|_2+\|f\|_{\frac65}\right)\|P\Delta w\|_2\,,
\ea
$$
and
$$\ba{rl}\medskip
\half\ode{}t\|\nabla w(t)\|_2^2+\|\partial_t w\|_2^2=&\!\!\!\!\!(\tilde{ v_{(\eta)}}\cdot\nabla w,\partial_t w)+(\tilde{f},\partial_tw)\\
\le&\!\!\!\!\!\left(\|\tilde{ v_{(\eta)}}\|_{L^{\infty,\infty}}\|\nabla w\|_2+\|f\|_{\frac65}\right)\|\partial_t w\|_2\,,\ea
$$
respectively. Therefore, by Cauchy--Schwartz inequality,  \eqref{2.7}, \eqref{2.6}, and the well-known estimate $\|w||_{2,2}\le c\,(\|P\Delta w\|_2+\|\nabla w\|_2)$ valid with a constant $c$ independent of $R$ \cite[Lemma 1]{Hey}, we show
\be
\int_0^t\big(\|\partial_\tau w(\tau)\|_2^2+\|w(\tau)\|^2_{2,2}\big)\le C\,,\ \ \mbox{$t\in [0,T]$, all $T>0$}\,, 
\eeq{2.8}
where the constant $C$ depends {\em only} on 
$T$, $
\|\tilde{ v_{(\eta)}}\|_{L^{\infty,\infty}}$, and $\|(v_0)^{(\eta)}\|_{1,2}\,,
$ 
and is therefore {\em independent of} $R$.Thus, coupling the classical Galerkin method together with the estimate \eqref{2.8}, we show the existence of a solution $(w_R,{\sf p}_R)$ to problem \eqref{2.5} in the class \eqref{2.2} (with $\real^3$ replaced by $B_R$). Note that this solution continues to satisfy the uniform bounds \eqref{2.7} and \eqref{2.8}. As a result, we may let $R\to\infty$ along a sequence and use \eqref{2.8}, to prove that $(w_R,{\sf p}_R)$ converges (in suitable topology) to the desired solution $(w_\eta,{\sf p}_\eta)$ for which \eqref{2.7} and \eqref{2.8} hold \cite[p. 660 and {\em ff}.]{Hey}.   Next, take $v_0\equiv0$. By H\"older inequality, \eqref{2.7} and \eqref{H}, we have
\be
\|\tilde{v_{(\eta)}}\cdot\nabla w_\eta\|_{L^{\frac43,\frac43}}\le \|\tilde{v_{(\eta)}}\|_{L^{4,4}}\|\nabla w_\eta\|_{L^{2,2}}\le c_3\,\|v\|_{L^{4,4}}\|f\|_{L^{2,\frac65}}\,,
\eeq{2.20}
which implies, in particular, 
$$
\tilde{v_{(\eta)}}\cdot\nabla w_\eta\in L^{\frac43,\frac43}\,.
$$
Therefore, from classical results (e.g. \cite[Theorem VIII.4.1]{GaBook}) the problem
\be\ba{cc}\smallskip
\partial_t w=\Delta w-\nabla {\chi}+F\,,\ \
\Div w=0\,
\ \ \mbox{in}\,\, \real^3\times (0,T)\,,
\\
w(x,0)=0\ \ x\in\real^3\,,
\ea
\eeq{2.21}
with $F:=\tilde{ v_{(\eta)}}\cdot\nabla w_\eta+\tilde f$ has at least one solution $(w_\star,\chi_\star)$ such that
$$ 
(w_\star,\nabla\chi_\star)\in \mathscr W^{1,\frac43}\times L^{\frac43,\frac43}\,.
$$
However, by uniqueness \cite[Lemma VIII.4.2]{GaBook} we must have $(w_\eta,\nabla{\sf p}_\eta)\equiv (w_\star, \nabla\chi_\star)$, which completes the proof of the lemma.\par\hfill$\square$.
\smallskip\par
{\bf\em Proof of Theorem 1.1}. Let $u_\eta$ be the solution of \lemmref{2.2} corresponding to $f\equiv 0$.
From \eqref{1.5}$_1$ and \eqref{2.1} we infer, for arbitrary $\varphi\in \cald_T$,
\be
\Int{0}{T}\big((v-u_\eta),\partial_t\varphi+\Delta\varphi+v_{(\eta)}\cdot\nabla\varphi\big)=\Int0T((v-v_{(\eta)})\cdot\nabla\varphi,v)-(v_0-(v_0)^{(\eta)},\varphi(0))\,.
\eeq{4.1}
Next, let $\Phi_\eta(x,t):=\psi_\eta(x,T-t)$, $\Xi_\eta(x,t):={\sf p}_\eta(x,T-t)$, $(x,t)\in\real\times [0,T]$, with $(\psi_\eta,{\sf p}_\eta)$  solution  constructed in \lemmref{2.2} corresponding to $v_0\equiv 0$. It then follows that $(\Phi_\eta,\nabla \Xi)\in\mathscr W^{1,\frac43}\cap \mathscr W^{1,2}\times L^{\frac42,\frac43}\cap L^{2,2}$,  and that $(\Phi,\Xi)$ solves the final-value problem 
\be 
\ba{cc}\smallskip
\partial_t\Phi+ v_{(\eta)}\cdot\nabla \Phi+\Delta \Phi=\nabla {\Xi}+f\,,\ \
\Div \Phi=0\,
\ \ \mbox{in}\, \real^3\times (0,T)\,,
\\
\Phi(x,T)=0\ \ x\in\real^3\,.
\ea
\eeq{4.2}
Since $u_\eta\in \mathscr W^{1,2}$, by embedding $u_\eta\in L^{4,4}$. Moreover, the trilinear form $\mathscr T:=\int_0^T(v_1\cdot\nabla v_2,v_3)$ is continuous in $L^{4,4}\times L^{2,2}\times L^{4,4}$. Thus, by  \eqref{H}, and with the help of a density result proved in Lemma A.1 in the Appendix, we can replace $\Phi_\eta$ for $\varphi$ in \eqref{4.1} and use \eqref{4.2}$_1$ to show
\be
\Int{0}{T}\big((v-u_\eta),f\big)=\Int0T((v-v_{(\eta)})\cdot\nabla\Phi_\eta,v)-(v_0-(v_0)^{(\eta)},\Phi_\eta(0))\,,
\eeq{4.3}
where we also observed that, by another  density argument based on \eqref{1.5}$_2$ and \eqref{2.1}$_2$ \cite[Theorem II.7.1]{GaBook}, it is 
$
\int_{0}^{T}\big((v-u_\eta),\nabla\Xi\big)=0\,.
$
We next pass to the limit $\eta\to 0$ in \eqref{4.3}. By \eqref{2.4},  $\|\Phi_\eta(0)\|_2$ is bounded by a constant independent of $\eta$. Therefore, since $v_0\in L^2_\sigma$,
\be
\lim_{\eta\to 0} (v_0-(v_0)^{(\eta)},\Phi_\eta(0))=0\,.
\eeq{4.4}
Next, by the continuity of the trilinear form $\mathscr T$,
 \eqref{H} and \eqref{2.4} it immediately follows that
\be
\lim_{\eta\to 0}\Int0T((v-v_{(\eta)})\cdot\nabla\Phi_\eta,v)=0\,.
\eeq{4.5}
Finally,  \lemmref{2.2}, and in particular \eqref{2.4}, entails the existence of an element $u\in L^{\infty,2}\cap L^2(0;T;H^1)$ such that, along a sequence $\{\eta_n\}$, 
$$ 
u_{\eta_n}\to u\,,\ \ \mbox{ $\textrm{weak}-\star$ in  $L^{\infty,2}$\, and\, weakly in  $L^{2}(0,T;H^1)$}\,.
$$
As a consequence, recalling that $f\in C_0^\infty(\real^3\times(0,T))$ we infer
\be
\lim_{n\to\infty}\Int{0}{T}\big((v-u_{\eta_n}),f\big)=\Int{0}{T}\big((v-u),f\big)\,.
\eeq{4.6}
Collecting \eqref{4.3}--\eqref{4.6} we thus conclude
$
\int_{0}^{T}\big((v-u),f\big)=0\,,
$
which, by the arbitrariness of $f$ completes the proof of the theorem.\par\hfill$\square$
\Br (a) By a slight modification of the proof just given, one can show the same result with the assumption \eqref{H} replaced by $v\in L^{\frac{2s}{s-3},s}$, $s>3$, thus reobtaining (by a different argument) a well known result of Fabes {\em et al.} \cite[Theorems (2.1) and (5.3)]{FJR}. In fact, by our method one could  cover also the borderline case $v\in C([0,T];L^3)$ that is excluded in \cite{FJR}.\par 
(b) The theorem continues to hold with  $\real^3$ replaced by a sufficiently smooth bounded or exterior domain, $\Omega$, and with $v$ vanishing at $\partial\Omega$ in the sense of  ``very weak" solutions; see, for the latter,  \cite{FarGaSoh,FarKoSoh}. In such a case, the proofs of \lemmref{2.2} (especially for $\Omega$ exterior) and Lemma A.1 in the Appendix may become technically (but not conceptually) more involved. This is also a reason why we preferred to treat here the Cauchy problem.\par  
(c) As an immediate corollary to \theoref{1}, we obtain the following {\em Liouville-type} result: If $v$ satisfies \eqref{H}, \eqref{1.5} with $v_0\equiv 0$, then $v\equiv 0$.  
\ER{1}
\section*{Appendix} 
\setcounter{equation}{0}
\renewcommand{\theequation}{A.\arabic{equation}}
{\bf Lemma A.1} {\sl Let $\mathscr W_0^{1,q}=\{\psi\in \mathscr W^{1,q}\cap C([0,T];H^1):\ \psi(T)=0\}$, $1<q<\infty$. Then $\cald_T$ is dense in $\mathscr W^{1,q}_{0}$.  }\smallskip\par
\noindent
{\em Proof.} Let $\chi_R\in C_0^\infty(\real^3)$ with $\chi_R(x)=1$ for $|x|\le R$,  $\chi_R(x)=0$ for $|x|\ge 2R$ such that
\be
|\nabla\chi_R(x)|\le C\,R^{-1}\,,\ \ |D^2\chi_R(x)|\le C\,R^{-2}\,,
\eeq{0}
with $C$ independent of $R$.
Likewise, for  small $h>0$,  let $\zeta=\zeta_h(t)$ be a $C^\infty$ function such that 
\be
\zeta_h(t)=\left\{\ba{ll}0\ \, \mbox{for $t\ge T-h$}\\
1\ \, \mbox{for $t\le T-2h$}
\ea\right. \,,\ \
|\zeta_h^\prime(t)|\le C\,h^{-1}
\eeq{1}
with $C>0$ independent of $h$. Set
$
\psi_{\eta_1,\eta_2}(x,t)=\int_0^Tj_{\eta_1}(t-s)\psi^{(\eta_2)}(x,s)ds\,,
$
and  denote by $w_R=w_R(x,t)$ a solution to the problem
\be
\Div w_R=-\nabla\chi_R(\zeta_\eta\psi)_{\eta_1,\eta_2}\,,\ \ w_R\in C_0^\infty(\real\times [0,T))\,,
\eeq{2}
satisfying the estimates
\be\ba{rl}\medskip
R^{-1}\|w_R\|_q+\|\nabla w_R\|_q\!\!\!&\le C\,\|\nabla\chi_R(\zeta_h\psi)_{\eta_1,\eta_2}\|_q
\\ \medskip
\|\nabla(\nabla w_R)\|_q\!\!\!&\le C\,\|\nabla(\nabla\chi_R(\zeta_h\psi)_{\eta_1,\eta_2})\|_q\,,\\ 
\|\partial_t w_R\|_q\!\!\!&\le C\,R\,\|\nabla\chi_R\,\partial_t(\zeta_h\psi)_{\eta_1,\eta_2}\|_q \,,
\ea
\eeq{3}
where the constant $C$ is independent of $R$. The existence of  a function $w_R$ is  known \cite[Theorem III.3.3 and Exercise III.3.6]{GaBook}. We shall  show that, for any $\varepsilon>0$ there are  sufficiently large $\bar R$ and small $\bar\eta$ such that, for all $R\ge \bar R$ and all $h,\eta_1,\eta_2\le\bar \eta$,
\be
\|w_R\|_{\mathscr W^{1,q}}+\max_{t\in[0,T]}\|w_R(t)\|_{1,2}<c_0\,\varepsilon\,.  
\eeq{4}
To this end, set $B{(R)}:=\{R<|x|<2R\}$. By \eqref{3}$_3$ and \eqref{0} we have
\be
\int_0^T\|\partial_tw_R\|_q^q\le C\int_0^T\big[\|(\zeta_h^\prime\psi)_{\eta_1,\eta_2}\|_{q,B(R)}^q+\|(\zeta_h\partial_t\psi)_{\eta_1,\eta_2}\|_{q,B(R)}^q\big]:=C\,(I_1+I_2)\,,
\eeq{5}
where the subscript $B(R)$ means that the spatial integration is restricted to $B(R)$. By the triangle inequality,  the properties of the mollifier and \eqref{1}$_1$, we get
\be
I_2\le \int_0^T\left[\|(\zeta_h\partial_t\psi)_{\eta_1,\eta_2}-\zeta_h\partial_t\psi\|_q^q+\|\partial_t\psi\|_{q,B(R)}^q\right]<\varepsilon\,,
\eeq{6}  
for all  sufficiently large $R$ and small $\eta_1,\eta_2$. Furthermore, also with the help of \eqref{1}, we show
$$
I_1\le \frac{c_1}h\int_{T-2h}^{T-h}\|\psi(\tau)\|_q^q\,.
$$
Thus, since $\mathscr W^{1,q}\subset C([0,T];L^q)$ and $\psi(T)=0$, this entails, with $h$ sufficiently small, 
\be
I_1\le {c_1}\,\varepsilon\,.
\eeq{7}
Next, by the properties  \eqref{0} of $\chi_R$,  for all sufficiently large $R$ and small $\eta_1,\eta_2$, it easily follows that
\be
\int_0^T\|w_R(t)\|_q^q \le C\,\big[\,\|(\zeta_h\psi)_{\eta_1,\eta_2}-\zeta_h\psi\|_{L^{q,q}}^q+\int_0^T\|\psi(t)\|_{q,B(R)}^q\big]<\varepsilon
\eeq{8}
In analogous (and simpler) fashion we can show
\be
\int_0^T\left(\|\nabla w_R(t)\|_q^q+\|D^2 w_R(t)\|_q^q\right)<\varepsilon
\,.
\eeq{9}
Finally, from \eqref{4}$_{1,2}$, \eqref{0} and again the triangle inequality we get
$$
\|w_R(t)\|_{1,2}\le C\,\left[\|(\zeta_h\psi)_{\eta_1,\eta_2}(t)-\zeta_h\psi(t)\|_{1,2}+\|\psi(t)\|_{1,2,B_R}\right]\,,\ \ t\in [0,T]\,,
$$
which, in turn, since $\psi\in C([0,T];H^1)$, for sufficiently large $R$ and small $\eta_1,\eta_2$, implies
$$
\max_{t\in[0,T]}\|w_R(t)\|_{1,2}<\varepsilon.
$$
As a result, \eqref{4} follows from the latter and \eqref{5}--\eqref{9}. Set
\be
\Psi_{R,h,\eta_1,\eta_2}:=\chi_R(\zeta_h\psi)_{\eta_1,\eta_2}-w_R\equiv U_{R,h,\eta_1,\eta_2}-w_R\,.
\eeq{10}
It is at once checked that $\Psi_{R,h,\eta_1,\eta_2}\in C_0^\infty(\real^3\times [0,T))$ and, in view of \eqref{2}, also that $\Div\Psi_{R,h,\eta_1,\eta_2}=0$. Consequently, $\Psi_{R,h,\eta_1,\eta_2}\in\cald_T$. In view of \eqref{10} and what we have already established in \eqref{4}, to complete the  proof, it remains to show that, for a given $\varepsilon>0$, it holds  
$$
\|U_{R,h,\eta_1,\eta_2}-\psi\|_{\mathscr W^{1,q}}+\max_{t\in [0,T]}\|U_{R,h,\eta_1,\eta_2}(t)-\psi(t)\|_{1,2}<\varepsilon\,,
$$
for all sufficiently large $R$ and small $h,\eta_1,\eta_2$. However, the proof of the latter is quite straightforward since, under the given assumptions on $\psi$, it only requires the use of classical properties of mollifiers and, therefore, it will be omitted.\par\hfill$\square$

\ed